\documentclass{amsart}
\usepackage{amsthm,amssymb,mathrsfs}
\newtheorem{theorem}{Theorem}

\newtheorem{definition}{Definition}

\newtheorem{example}{Example}
\begin{document}
\title[Multiplicity of zeros]{ Multiplicity of zeros and discrete orthogonal polynomials}
    
\author[I. Krasikov]{Ilia Krasikov}
     
\address{   Department of Mathematical Sciences,
            Brunel University,
            Uxbridge UB8 3PH United Kingdom}
\email{Ilia.Krasikov@brunel.ac.uk}
			      
			       
\newcommand{\majp}[1]{\succ\!\!_{#1}}
\newcommand{\ee}{\end{equation}}
\newtheorem{remark}{Remark}
\newtheorem{cor}{Corollary}
\newtheorem{conjecture}{Conjecture}
\newcommand{\cD}{\mathcal{D}}
\newcommand{\cS}{\mathcal{S}}
\newcommand{\cP}{\mathcal{P}}
\newcommand{\cG}{\mathcal{G}}
\newcommand{\ZN}{\mathbb{Z}}
				
\begin{abstract}
We consider a problem of bounding the maximal possible multiplicity of a zero
of some expansions $\sum a_i F_i(x),$
at a certain point $c,$ depending on the chosen family $\{F_i \}.$
The most important example is a polynomial with $c=1.$
It is shown that 
this question
naturally leads to
discrete orthogonal polynomials. 
Using this connection we
derive some new  bounds, in particular  on the 
multiplicity of the zero at one of a polynomial
with a prescribed norm.
\end{abstract}

\maketitle
\section{Introduction}
This paper deals with an old
question about the maximal
possible multiplicity $\mu$ of a zero which a polynomial with
a prescribed norm can have at 1.
The first important results were obtained by Bloch, Polya and Schur in early 30-s \cite{bp}, \cite{schur}.
Recently Peter Borwein and  T\'{a}mas Erd\'{e}lyi revived interest in the problem
and substantially improved our understanding of the matter, 
see \cite{bor1}, \cite{bor}, \cite{erdel1} and the references therein.
Although the most interesting problems arise when the coefficients
of the polynomial are integers restricted to a few values, e.g. $\{-1,0,1\},$
today we hardly know how to exploit those constraints apart from a few cases, such as the result of Boyd \cite{boyd}.
Polynomials with a high vanishing at 1 naturally appear
in coding theory. Since this seemingly important connection belongs rather to folklore, let us describe it
for the simplest binary case. The relevant definitions can be found in \cite{MWS}.
Let $C$ be a binary code of length $n$ and size $|C|,$
with the distance and the dual distance distributions $1=B_0, ...,B_n,$
and $B'_0,...,B'_n,$ respectively. They are related by the McWilliams transform
$$|C| \, B'_i = \sum_{i=0}^n B_j K_i^n (j),$$
where $K_i^n (j),$ are the binary Krawtchouk polynomials, defined by the generating function
$$(1-z)^x (1+z)^{n-x}=\sum_{j=0}^{\infty} K_j^n(x)z^j.$$
Multiplying both sides by $x^i$ and summing up, one gets
$$|C| \sum_{i=0}^n B'_i x^i=\sum_{i=0}^n B_i (1-x)^i (1+x)^{n-i}.$$
If the code distance is $d$, that is $B_1=B_2=...=B_{d-1}=0,$
then 
$$\sum_{i=0}^n \left( |C|B'_i - {n \choose i} \right) =(1-x)^d \; \sum_{i=d}^n B_i (1-x)^{i-d} (1+x)^{n-i},$$
yielding the sought polynomial. Unfortunately, the distance distribution is usually
unknown, and applications of this idea are rather tricky, we refer to \cite{abl} and the references therein.

In this paper we describe how using discrete orthogonality one can derive bounds
on $\mu$ at a certain point $c$ for some more general expansions $\sum a_i F_i$.
In fact, the point $c$ will depend on the chosen basis $\{ F_i \},$
e.g. $c=1$ for $F_i=x^i$ and 0 for $F_i=L_i^{(\alpha )}(x),$ the Laguerre polynomials.
On the other hand, the obtained bounds are independent of the basis
and are the same as for polynomials. The definition of the feasible functions $F_i$
is rather technical and will be given later on.
Our main result is Theorem \ref{theogl}. As an easy corollary we establish
some explicit inequalities on $\mu$ in weighted
$\ell_2$ and $\ell_{\infty}$ norm.

Throughout the paper we will use the following notation.
Let ${\cP}^{m}$ be the set of all complex polynomials of degree $m$ or less.
We will use $\ZN_{m,n}$ for the set of integers $m,m+1,...,n.$
The multiplicity of the zero of $f(x)$ at $c  \in \mathbb {C}$
is denote by $ \mu (f,c).$ Background information on classical discrete orthogonal polynomials
can be found in  \cite{szego} and \cite{niki85}.
\section{Discrete Orthogonality}
We start with the following definition describing a feasible system of functions
which can be tackled by our method.
\begin{definition}
\label{def1}
Let $F_0(x), F_1(x),...F_{n}(x),$ be a sequence of functions analytic in
an open disc centered at $c.$
It will be called a $\Delta(n, c) $ basis 
if there are $n+1$ polynomials $R_0 (x), R_1(x),...,R_{n}(x),$ such that\\
(i)   the derivatives $F_i^{(j)}(c)=R_j(i),$ $i,j=0,1,..,n,$\\ 
(ii)  $R_i (x)$ form a basis in  ${\cP}^{m}.$
\end{definition}
Let us consider a few examples.
\begin{example}
\label{ex1}
$F_i(x)=x^i$ are a $\Delta(n,1)$
basis for any $n \in \ZN^+$, with 
$$R_j(x)=\Gamma (x+1)/\Gamma (x-j+1).$$
\end{example}
\begin{example}
\label{ex2}
$F_i(x)=(1-x)^i (1+x)^{(n-i)}$ are a $\Delta(n, 0) $ basis 
with $R_j(x)=j! K_j^n (x),$ where $K_j^n (x)$ are the binary Krawtchouk
polynomials.
\end{example}
\begin{example}
\label{ex3}
$$F_i(x)=\frac{L_i^{(\alpha )}(x)}{L_i^{(\alpha )}(0)}, \; \; \alpha \ne -1,-2,...-n,$$
where $L_i^{(\alpha )}(x)$ are  the Laguerre polynomials,
are a $\Delta( n,0)$ basis for any $n \in \ZN^+$, with 
$$R_j(x)=(-1)^j {x \choose j}{j+ \alpha \choose j}^{-1}.$$
This readily follows from the well known properties of the Laguerre
polynomials \cite{szego}:
$$\frac{d}{dx}L_i^{(\alpha )}(x)=-L_{i-1}^{(\alpha +1 )}(x), \; \; L_i^{(\alpha )}(0)={i+ \alpha \choose i}.$$
\end{example}
\begin{example}
\label{ex4}
Generating functions
$$\phi (x,z)=\sum_{i=0} \frac{p_j(x)}{j!} (z-c)^j,$$
for different families of polynomials $p_j,$ such that $deg(p_j)=j,$
give many examples of,
not necessarily polynomial, 
$\Delta(n, c) $ bases
$\{ \phi (i,c)\}, \; i=0,1,...,n,$
with $R_j(x)=p_j(x).$
\end{example}
Now we should introduce some definitions and notation related to
discrete orthogonal polynomials.
Denote by $\cG (w,g,\cD )$, 
a family of discrete orthonormal polynomials
$g_0(x), g_1(x),...;$  $deg(g_i)=i,$ on
$\cD  \subseteq \ZN ,$
with respect to the weight function $w(x),$ $w(x) >0$ iff $x \in \cD .$ 
That is
$$ \sum_{x \in \cD } w(x)g_i(x)g_j(x)=\delta_{ij}.$$
Let $\{ F_i(x) \}$ be a $\Delta(n,c)$ basis.
For given $\cG (w,g,\cD )$ and $p(x)=\sum_{i=0}^n a_i F_i(x),$ 
we set $\cS = \{i \, | \, a_i \ne 0 \},$
and write $w \succ p$ if $ \cS \subseteq \cD .$
We also denote by $A_w (x)$ 
the polynomial of the least degree such that $\frac{a_i}{w(i)} =A_w (i)$ for $i \in \cD , $
and put $N=deg (A_w (x)).$ 
Thus, $A_w (x)$ is just the corresponding interpolation polynomial.
Notice that in general,
$N$ can be greater then $n$ since one has to put $ a_i=0$ for $i \in \cD \smallsetminus \cS .$
\begin{theorem}
\label{thgl}
Let $p(x)=\sum_{i=0}^n a_i F_i(x),$ where $\{F_i(x) \}$ is a $\Delta(n,c)$ basis. Then
$\mu (p,c)= \mu \ge 1 ,$ iff either of the following conditions hold\\
(i) for any $f(x) \in {\cP}^{\mu-1},$
\begin{equation}
\label{glav}
\sum_{i =0}^n a_i f(i)=0.
\end{equation}
(ii) for $j \in \cD,$
\begin{equation}
\label{zaglav}
\frac{a_j}{w(j)}=\sum_{k= \mu}^N \lambda_k g_k (j)
\end{equation}
where $g_i \in \cG (w,g,\cD )$ and $w \succ p.$
\end{theorem}
\begin{proof}
The condition $\mu (p,c)= \mu \ge 1 ,$ means $p(c)=p'(c)=...=p^{(\mu -1)}(c)=0.$
Our definition of a  $\Delta(n,c)$ basis implies  $p^{(j)}(c)=\sum_{i=0}^n a_i R_j(i),$
and, since $f(i)$ has a unique representation 
as a linear combination of $R_j(i), \; j=0,...,\mu -1,$
we get for appropriate $\lambda_j ,$
$\sum_{i =0}^n a_i f(i)=\sum \lambda_j p^{(j)}(c)=0.$
The second claim follows from (\ref{glav}) and the expansion $\frac{a_j}{w(j)}=A_w(j)=\sum_{k=0}^N \lambda_k g_k (j),$ 
by the orthogonality.
\end{proof}
\begin{theorem}
\label{theogl}
Let $p(x)=\sum_{i=0}^n a_i F_i(x),$ where $\{F_i(x) \}$ is a $\Delta(n,c)$ basis, and
let $\mu (p,c)= \mu \ge 1 .$ 
Then\\ 
(i)  for any $\cG (w,g,\cD )$ such that $w \succ p,$ and $s \in \cD ,$
the following sharp inequality holds
\begin{equation}
\label{ozl2}
w^2(s) \, \sum_{i \in \cD } \frac{|a_i|^2}{w(i)} \ge |a_s|^2 \, \left( \sum_{j= \mu}^N g_j^2(s) \right)^{-1}
\end{equation}
(ii) for any $\cG (w,g,\cD )$ such that $w \succ (p(x)- a_s F_s(x)),$ and $s \notin \cD ,$
\begin{equation}
\label{condg2}
\max_{k \in \cD} \frac{|a_{k}|}{w(k)} \ge
|a_s| \sum_{j=0}^{\lfloor \frac{\mu-1}{2} \rfloor } g_j^2(s)
\end{equation}
\end{theorem}
\begin{proof}
(i) Putting $A_w(i)=\sum_{j= \mu }^N \lambda_j g_j(i),$
we get
$$\sum_{i \in \cD } \frac{|a_i|^2}{w(i)}
=\sum_{x \in \cD} w(x) \, |A_w (x)|^2
=\sum_{x \in \cD } w(x) | \sum_{i=\mu}^N \lambda_i g_i(x) |^2
=\sum_{i=\mu}^N | \lambda_i |^2 $$
We also have
\begin{equation}
\label{usl}
a_s=w(s) \sum_{i=\mu}^N \lambda_i g_i(s)
\end{equation}
The minimum of the quadratic form $\sum_{i=\mu}^N | \lambda_i |^2 $ subjected to (\ref{usl})
is
$$\frac{|a_s |^2}{w^2(s)} \left( \sum_{j= \mu}^N g_j^2(s) \right)^{-1},$$
yielding (\ref{ozl2}).\\
(ii) This follows from a well-known extremal property of the Christoffel function. Namely,
let $r=\lfloor \frac{\mu-1}{2} \rfloor ,$ $M=\max_{k \in \cD} \frac{|a_{k}|}{w(k)},$ and $q(x)= \sum_{j=0}^r g_j(s)g_j(x).$ 
Putting $f(x)=q^2(x)$ in (\ref{glav}), we get
$$|a_s | \, q^2(s)  \le \sum_{ k \ne s} |a_k | \, q^2(k),$$
hence
$$
|a_s| \, q^2(s)  \le  \sum_{ k \ne s} |a_k | \, q^2(k)=
\sum_{ k \ne s}  \frac{|a_{k}|}{w(k)} \, w(k)\, q^2(k) \le 
M \, \sum_{k \in \cD} w(k)\, q^2(k)=  $$
$$
M \sum_{k \in \cD} \sum_{j,l=0}^r w(k)g_j(k)g_l(k)g_j(s)g_l(s)=M \sum_{j,l=0}^r g_j(s)g_l(s)=M \, q(s).
$$
Thus $M \ge |a_s | q(s),$ what is equivalent to (\ref{condg2}).
\end{proof}
\section{Corollaries}
There are four families of classical discrete orthogonal polynomials which can be readily used in
Theorem \ref{theogl}. The first two are the Hahn and the Krawtchouk polynomials both orthogonal on $\ZN_{0,n}.$
Let us write down the corresponding formulae for the orthonormal case \cite{niki85}.\\
{\bf Hahn Polynomials} $Q_k^n (x, \alpha , \beta)$,
$\alpha , \beta >-1$ or  $\alpha , \beta <-n,$
$$w(x)= {x+\alpha \choose x}{n-x+ \beta \choose n-x}.$$
\begin{equation}
\label{hahn}
Q_k^n (x, \alpha , \beta)= \sqrt{d_k}\, \sum_{j=0}^k (-1)^{j}
\frac{{k \choose j}{k+\alpha+\beta+j \choose j}{x \choose j}}{{j+\alpha \choose j}{n \choose j}} \, ,
\end{equation}
where
$$d_k=\frac{ (2k+ \alpha + \beta+1){k+ \alpha \choose k}{n \choose k}}{(n+1) {k+ \beta \choose k}{n+k+\alpha +\beta +1 \choose n+1}}$$
\begin{equation}
\label{hahn0}
\left(Q_k^n (0, \alpha , \beta) \right)^2 = d_k
\end{equation}
An important special case is the
discrete Chebyshev polynomials
$ T_k^n(x)=Q_k^n(x,0,0),$
corresponding to $w(x)=1.$ 
\begin{equation}
\label{dcheb0}
\left(T_k^n(0)\right)^2 =\frac{(2k+1)n!^2 }{(n-k)!(n+k+1)!}
\end{equation}
\begin{equation}
\label{dcheb1}
\left(T_k^n(-1)\right)^2 =\frac{(2k+1)(n-k)!(n+k+1)!}{(n+1)!^2}
\end{equation}
{\bf Krawtchouk Polynomials} $K_k^n(x,q)$, $0 <q < \infty, \; \; w(x)= {n \choose x}q^x .$
\begin{equation}
\label{kraw}
K_k^n(x,q)= \sqrt{{n \choose k}^{-1} q^{k} (1+q)^{-n}} \; \sum_{j=0}^k (-q)^{-j} {x \choose j}{n-x \choose k-j},
\end{equation}
\begin{equation}
\label{kraw0}
\left( K_k^n (0,q) \right)^2={n \choose k}q^k (1+q)^{-n}
\end{equation}
Although the substitution of the above expressions in (\ref{ozl2}) and (\ref{condg2}) gives formally 
the sought result, it has a little meaning without the corresponding asymptotic.
To avoid cumbersome formulae we will consider only
a few cases which allow relatively simple answer. We formulate the results in terms of $a_0$,
of course one can easily restate them in terms of $a_n$.\\
We also make use of the following elementary inequality, which can be readily proved
by induction on $k.$
\begin{equation}
\label{oze}
\frac{(n-k)!(n+k)!}{n!^2} \ge e^{2k^2/(2n+1)}, \; \; 0 \le k \le n.
\end{equation}
\begin{theorem}
\label{thvtor}
Let $p(x)=\sum_{i=0}^n a_i F_i(x),$ where $\{F_i(x) \}$ is a $\Delta(n, c) $ basis, and let
$\mu = \mu (p,c) \ge 1 .$ Then
\begin{equation}
\label{eq1}
\sum_{i=0}^n |a_i |^2 \ge \frac{(n-\mu )!(n+\mu )!}{n!^2} |a_0|^2 \ge 
e^{2 \mu^2/(2n+1)}|a_0|^2 
\end{equation}
\begin{equation}
\label{eq2}
1+\max_{k \ge 1}|  \frac{a_k }{a_0}|  \ge
\frac{(n- \lfloor \frac{\mu +1}{2} \rfloor )!(n + \lfloor \frac{\mu +1}{2} \rfloor )!}{n!^2} 
\ge \exp{\left( \frac{2 {\lfloor \frac{\mu +1}{2} \rfloor}^2}{2n+1}\right)}
\end{equation}
\begin{equation}
\label{eq3}
\sum_{i=0}^n \frac{|a_i |^2 q^{-i}}{{n \choose i}}
\ge \frac{|a_0 |^2}{ \left(1-(1+q)^{-n} \sum_{i=0}^{\mu -1} {n \choose i}q^i \right)}
\end{equation}
for $0 < q <\infty .$
\end{theorem}
\begin{proof}
First notice that the last two inequalities in (\ref{eq1}) and (\ref{eq2}) follow from (\ref{oze}).\\
To prove (\ref{eq1}) we choose  $g_i(x)=T_i^n(x)$ in (\ref{ozl2}). 
Using (\ref{dcheb0}), we get, with the convention $(-1)!=0$,
$$\sum_{i=\mu}^n g_i (0)^2 = \sum_{i=\mu}^n \frac{(2i+1)n!^2}{(n-i)!(n+i+1)!}=$$
$$ \sum_{i=\mu}^n \left(\frac{n!^2}{(n-i)!(n+i)!} -\frac{n!^2}{(n-i-1)!(n+i+1)!} \right)=
\frac{n!^2}{(n-\mu )!(n+ \mu )!}.$$
By (\ref{ozl2}) this yields (\ref{eq1}).
To prove (\ref{eq2}) we consider $g_k(x) =T_k^{n-1} (x-1).$
This system is orthogonal on $\ZN_{1,n}$ and
(\ref{dcheb1}) yields
$$n!^2 \sum_{i=0}^r (g_i (0))^2= \sum_{i=0}^r (2i+1)(n-i-1)!(n+i)!=$$
$$ \sum_{i=0}^r \left( (n-i-1)!(n+i+1)! -(n-i)!(n+i)! \right)
=(n-r-1)!(n+r+1)! -n!^2.  $$
Choosing $r=\lfloor \frac{\mu -1}{2} \rfloor ,$ we get (\ref{eq2}) by (\ref{condg2}).
Using Krawtchouk polynomials and (\ref{kraw0}) we obtain (\ref{eq3}) by (\ref{ozl2}).
\end{proof}
Notice that(\ref{eq2})
gives the best currently known bound on $\mu$ for polynomials with coefficients $\{-1,0,1 \},$ \cite{fk}.
Namely, we have $\mu \lesssim 2 \sqrt{ n \ln{2}}.$

Two other families of discrete orthogonal polynomial are the Meixner and the Charlier polynomials
orthogonal on $\ZN_{0, \infty} .$ 
It is easy to check that the interpolating function $A_w (x)$, which is not a polynomial now, converges
as only finitely many $a_i \ne 0.$ Thus the only difference with the previous case is
that our bounds cannot be attained by polynomials and all inequalities are strict.
Again, we start with listing the relevant formulae.\\
{\bf Meixner Polynomials} $M_k(x, \beta, q),$ $\beta >0$ and $0 <q <1,$
$w(x)={x+\beta-1 \choose x} q^x .$
\begin{equation}
M_k(x, \beta, q)=(-1)^k \sqrt{\frac{(1-q)^{\beta} q^k}{{k+\beta-1 \choose k}}} \, 
\sum_{j=0}^k {x \choose j}{-x- \beta \choose k-j}q^{-j}
\end{equation}
\begin{equation}
\label{m0}
\left(  M_k(0, 1, q)\right)^2 =(1-q) q^k 
\end{equation}
\begin{equation}
\label{m1}
\left( M_k(-1, 1 ,q) \right)^2 = (1-q)q^{-k} 
\end{equation}
{\bf Charlier Polynomials} $C_k(x, \lambda )$, $ \lambda  >0,$ $w(x)=\frac{\lambda^x}{x!}.$
\begin{equation}
C_k(x,  \lambda )= \sqrt{\frac{\lambda^k e^{-\lambda}}{k!}} \, \sum_{i=0}^k (-\lambda )^{-i}{k \choose i}{x \choose i} i! 
\end{equation}
\begin{equation}
\label{c0}
\left( C_k(0, \lambda ) \right)^2 = \frac{\lambda^{k}e^{-\lambda}}{ k!}.
\end{equation}
Using (\ref{m0}),(\ref{m1}) and (\ref{c0}), similarly to the previous theorem we obtain
\begin{theorem}
Let $\{F_i(x) \}$ be a $\Delta(n, c) $ basis, 
$p(x)=\sum_{i=0}^{n} a_i F_i(x),$ 
and $\mu = \mu (p,c) \ge 1 .$ Then 
\begin{equation}
\sum_{i=0}^{n} |a_i|^2 q^{i} > |a_0 |^2 q^{\mu} , \; \; \; q >1
\end{equation}
\begin{equation}
\max_i { |a_i | q^{i}} >  \left( q^{ \lfloor \frac{\mu -1}{2} \rfloor } -\frac{1}{q}  \right) |a_0|,   \; \; \; q >1
\end{equation}
\begin{equation}
\sum_{i=0}^{n} |a_i|^2 q^{-i} i! > \left( 1-e^{-q} \, \sum_{i=0}^{\mu-1} \frac{q^i}{i!} \right)^{-1} |a_0|^2 ,  \; \; \; q >0
\end{equation}
\end{theorem}
It is interesting to compare the obtained bounds, especially (\ref{eq1})
and (\ref{eq2}), with the following classical result of Schur \cite{schur} (see
also \cite{erd} and \cite{bor1}). 
Let $\nu$ be the number of real roots of $p(x)=\sum_{i=0}^m a_i x^i ,$ then
$$ \sum_{i=0}^n |a_i |^2 \ge 2 a_0 a_n \, e^{(\nu^2- \nu)/2n}, $$
$$ \sum_{i=0}^n |a_i | \ge \sqrt{a_0 a_n } \; e^{\nu^2 /4n}.$$
Although we are not aware of any general claim of this type,
one could expect that upper bounds on both numbers $\nu$ and $\mu$ should be very close
for a reasonable choice of the norm.\\
It worth also noticing that one can view Theorem \ref{thgl} as giving
bounds on  discrete orthogonal polynomials. It would be very important
to learn how to construct polynomials with a small norm and a large
value of $\mu.$

\end{document}